\documentclass[12pt]{article}

\usepackage{amsmath,amsthm,amsfonts,amssymb}
\evensidemargin
.25in \textheight 8.5in \topmargin 0in \title{Solvable 
subgroups of
$Out(F_n)$ are virtually abelian} \author{Mladen 
Bestvina, Mark
Feighn, and Michael Handel\thanks{All three 
authors gratefully
acknowledge the support of the National Science 
Foundation.}}
\date{June 1997} 

\newtheorem{thm}{Theorem}[section] 
\newtheorem{lemma}[thm]{Lemma}
\newtheorem{cor}[thm]{Corollary}
\newtheorem{proposition}[thm]{Proposition}

\theoremstyle{definition}
\newtheorem{defn}[thm]{Definition}

\newtheorem{ex}[thm]{Example}

\theoremstyle{remark}

\begin{document}
\maketitle
\def\Cal{\cal}

\newcommand{\h}{\frak H}
\newcommand{\x}{{\Cal X}}
\newcommand{\ve}{{\Cal V}}
\newcommand{\PG}{PG_0(1)}
\newcommand{\vs}{\x}
\newcommand{\tbar}{\overline T}
\newcommand{\e}{\epsilon}

\def\G{{\cal H}}
\def\H{{\cal H}_0}
\def\R{\cal R}
\def\f{{F_n}}
\def\U{UPG_{\f}}
\def\Beta{\cal B}

\newcommand{\n}{n}
\newcommand{\bd}{\partial}

\def\hom(#1){H_1(#1,\Bbb Z/3\Bbb Z)}
\def\l(#1,#2){\ell_{#1}(#2)}
\def\roster{\begin{enumerate}}
\def\endroster{\end{enumerate}}
\def\definition{\begin{defn}}
\def\enddefinition{\end{defn}}
\def\subhead{\subsection\{}
\def\endsubhead{\}}
\def\head{\section\{}
\def\endhead{\}}
\def\example{\begin{ex}}
\def\endexample{\end{ex}}
\def\ves{\vs}
\def\cv{CV}
\def\ttt{G}
\def\ttm{f}
\def\oa{{\cal O}}
\def\o(#1){Out(#1)}
\def\F(#1,#2){{\cal F}_{#1}(#2)}
\def\Z{{\mathbb Z}}
\def\A{{\cal A}}
\def\oaf{\oa}
\def\oag{\oaf'}
\def\oah{\oaf'}
\def\ttmphi{\ttm}
\def\u{u}
\def\a{\mu}
\def\P{\Pssi}\def\Pssi{\Psi}
\def\g{\oag}
\def\tmm{f}
\def\glthree{GL(n,\Z/3\Z)}
\def\foa{f_\oaf}

\def\ti{\tilde}
\def\kr{upper triangular representative}
\def\pssi{\oa}
\def\phhi{\oa_0}
\def\Phhi{\Psi}
\def\h'{\H'}
\def\Q{Q}
\def\fG{f : G \to G}
\def\k{{\cal K}}
\def\akr{abelian Kolchin representative}
\def\he{FHE(G,{\cal V})}
\def\h{{\cal H}}
\def\egs{exponentially growing stratum}
\def\c{C_{\infty}}

\section{Introduction} \label{intro}

Let $\f$ denote a
free group of rank $n$. The group $\o(\f)$ contains
mapping class groups of
 compact
surfaces  and maps onto $GL(\n,\Z)$. It is perhaps
 not surprising
that $\o(\f)$ behaves at times like a mapping class
group and at times
like a linear group. J. Birman, A. Lubotzky, and J.
McCarthy
\cite{blm:s2a} showed that solvable subgroups of
mapping class groups
are finitely generated and virtually abelian. Of course,
$GL(3,\Z)$
contains the Heisenberg group which is solvable but
not virtually abelian. In
 this paper, we show
that, with respect to the nature of solvable subgroups,
$\o(\f)$ behaves more
like mapping class groups.

\begin{thm}\label{main} Every solvable subgroup of
$\o(\f)$
has a subgroup of index at most $3^{5n^2}$ that is
finitely generated and free
abelian.\end{thm}\noindent The rank of an abelian
subgroup
of $\o(\f)$ is bounded by $vcd(\o(\f))=2\n-3$ for
$\n>1$ \cite{cv:moduli}.

Since $Aut(\f)$ embeds in $Out(F_{n+1})$, solvable
subgroups of
$Aut(\f)$ are also virtually abelian. 
Theorem~\ref{main} complements
\cite{bfh:tits2} where we show that $\o(\f)$ satisfies
the Tits
Alternative, i.e. that subgroups of $\o(\f)$ are either
virtually solvable
or contain a free group of rank 2.

H. Bass and A. Lubotzky \cite{bl:niltech} have shown
that
solvable subgroups of $Out({\f})$ are virtually
polycyclic. In particular, they
 are finitely
generated.   We include an independent proof of this
fact for completeness and because
the ingredients of our proof are needed for the proof
of Theorem~\ref{main}.

        The starting point for this paper  is a short exact
sequence from \cite{bfh:tits1} and \cite{bfh:tits2}.
Begin with a solvable subgroup
$\h$ of $\o(\f)$.  After passing to a finite index
subgroup we may assume that  $\h$
acts trivially on
$H_1(\f;\Z/3\Z)$.  By Theorem 8.1 
of
\cite{bfh:tits1} and Proposition 3.5 of
\cite{bfh:tits2}   there is  an exact sequence
$$1\to\H\to\G\overset
\Omega \to\Z^b\to 1$$ where  $\H$ is $UPG$
(defined below).

     There are two parts to the proof of
Theorem~\ref{main}.  First we show that $\H$ is
abelian by constructing an embedding $\Phi : \H \to
{\mathbb Z}^r$.   Then we show that $\Phi$ extends 
to a
homomorphism $\Phi' : \h \to {\mathbb Z}^r$.
The
direct sum of $\Omega$ and $\Phi'$ is  an
embedding of $\h$ into ${\mathbb
Z}^{b+r}$ showing that $\h$ is finitely generated and 
free abelian.

     Our approach is motivated by  the special case that
$\h$ is realized as a subgroup of the mapping class
group of a compact surface $S$.  The
surface $S$ decomposes into a union of
annuli $A_1,\dots,A_r$ and subsurfaces $S_i$ of
negative Euler characteristic;
virtually every $\eta \in
\h$ is represented by a homeomorphism $f : S \to S$
that restricts to a Dehn twist on
each $A_j$ and that preserves each $S_i$.  If $\eta
\in \H$ then each $f|S_i$ is the
identity. The homomorphisms $\Phi$ and $\Phi'$  
are defined
by taking their $j^{th}$ coordinates to be   the
number of twists that occurs across  $A_j$. For a
further discussion of the
geometric case, see Example~\ref{Dehn}.

    We study an outer automorphism   ${\cal \eta}$
through its  lifts $\hat
\eta : \c \to \c$ to the Cantor set  at infinity $\c$ or
equivalently (see
Subsection~\ref{lifts}) through the automorphisms
$\Phi : F_n \to F_n$ that
represent $\eta$.  The $\c$ point of view simplifies
certain
proofs because it allows us to consider
 fixed
\lq directions\rq\  in $F_n$ that are not periodic
and therefore do not come from
fixed elements of $F_n$.

    In the course of proving Theorem~\ref{main} we
prove the following result which is
of independent interest.  We state it here in terms of
automorphisms although we prove
it in terms of  $\c$.

\begin{proposition}\label{lifting to Aut}  Every
abelian subgroup
$\h \subset Out(F_n)$ has a virtual lift $\ti \h
\subset Aut(F_n)$.  If $\gamma$ is
a non-trivial primitive element of $F_n$ that is
fixed, up to conjugacy, by each
element of $\h$  then $\ti
\h$ can be chosen so that each element of $\ti \h$
fixes $\gamma$.
\end{proposition}

   The paper is organized as follows.  In
section~\ref{prelim} we establish notation and
record known results for future reference.  In
sections~\ref{A} and \ref{fin gen} we
prove that
$\H$ is finitely generated and free abelian and that 
the
 above exact sequence is
virtually central.  In
section~\ref{lifting} we prove
Proposition~\ref{lifting to Aut} and in
section~\ref{proof of main} we prove
Theorem~\ref{main}.

\section{Notation and Preliminaries} \label{prelim}

\subsection{Lifts to $\c$}  \label{lifts}

   We assume that $F_n$ is identified with
$\pi_1(R_n,*)$ where $R_n$ is the rose with $n$
petals and with vertex
$*$.  Let  $\tilde R_n$  be the universal cover of
$R_n$ and let $\ti
*$ be a preferred lift of $*$.  The space of ends of
$\tilde R_n$ is a Cantor set that we denote $\c$.

 A {\it marked graph} is a graph $G$ with a preferred
vertex $v$
along with a homotopy equivalence $\tau : (R_n,*)
\to
(G,v)$ that identifies $\pi_1(G,v)$ with
$\pi_1(R_n,*)$ and so also
with $F_n$.

  Denote the space of ends of
$\Gamma$ by ${\cal E}(\Gamma)$.   The
marking homotopy equivalence lifts to an
equivariant map $\tilde \tau :
(\tilde R_n,\ti *) \to (\Gamma,\ti v)$  that induces
a
homeomorphism from $\c$ to ${\cal E}(\Gamma)$
(See for example section 3.2 of
\cite{bfh:tits1}). We use this homeomorphism to
identify  ${\cal E}(\Gamma)$ with $\c$
and so for the rest of this paper refer to the space of
ends of $\Gamma$ as $\c$.

   An outer automorphism $\eta$
of $F_n$ can be represented, in many ways, by a 
homotopy equivalence $f : G \to G$ of a marked 
graph.
 More precisely, $f : G \to G$ can be chosen so that 
when 
 $\pi_1(G,v)$ is identified with $F_n$, the 
outer automorphism of
$\pi_1(G,v)$ determined by $f$ agrees with $\eta$. 

 A preferred vertex $\ti v $ in the
universal cover
$\Gamma$  of $G$ provides an identification of the
group ${\cal T}$ of  covering
translations  of $\Gamma$ with $\pi_1(G,v)$ and so
with $F_n$.   The
action of
${\cal T}$ on $\Gamma$ extends to an action (of
$F_n$) on $\c$ by homeomorphisms. If $T$
is non-trivial, then the endpoints of its axis  are the
only fixed points for the action
of
$T$ on $\c$.

   Suppose that $f : G \to G$ represents $\eta$ and that 
$\ti f  :
\Gamma \to  \Gamma$ is a lift of $f$. For each $T
\in {\cal T}$ there exists a unique $T' \in
{\cal T}$ such that $\ti f T = T'
\ti f$. This defines an  automophism  $T \mapsto
T'$ of
${\cal T}$ and so an  automorphism $\Phi : F_n
\to F_n$.  It is easy to check that  the outer 
automorphism
class of $\Phi$ is $\eta$ and that as
$\ti f$ varies over all lifts of $f$, $\Phi$ varies
over all automorphisms representing
$\eta$. The subgroup
${\cal T}(\ti f)
\subset {\cal T}$ of covering translations that
commute with $\ti f$ corresponds to the
fixed subgroup
$Fix(\Phi) \subset F_n$.  Since $Fix(\Phi)$ is
quasiconvex and finitely
generated,  the closure in
$\c$ of the endpoints of axes of elements of
${\cal T}(\ti f)$ is  identified with the space of ends of
$Fix(\Phi)$ \cite{co:bcc}.

 Each  $\ti f : \Gamma \to \Gamma$ extends to a
homeomorphism $\hat f: \c \to \c$ . Denote the  
group of $F_n$-equivariant homeomorphisms of 
$\c$ by $EH(\c)$. The composite 
$\Phi \mapsto \ti f \mapsto \hat f$ defines an 
injective homomorphism
from $Aut(F_n)$ to  $EH(\c)$ that is independent of 
the choice of homotopy
equivalence $f : G \to G$ representing $\eta$. (See for 
example section 3.2 of
\cite{bfh:tits1}).  We will
sometimes write
$\hat
\eta$ instead of
$\hat f$ where $\eta
\in Out(F_n)$ is the outer automorphism
determined by $f : G \to G$.

  If
$\ti L \subset \Gamma$ is a line with endpoints $P$ 
and $Q$, we denote the line
with endpoints $\hat f(P)$ and $\hat f(Q)$ by $\ti 
f_\#(\ti L)$. If $\ti L$ is
the axis of $T \in {\cal T}$ then  $\ti f_\#(\ti L)$ is 
the axis of the
covering translation $T'
\in {\cal T}$ satisfying $\ti f T = T'
\ti f$.

  We conclude this subsection by recording some facts
for future reference. If $\ti \h
\subset Aut(F_n)$ is a lift of $\h \subset Out(F_n)$
then we denote the corresponding
lift to $EH(\c)$ by $\hat \h$. For any
subgroup
${\mathbb F}
\subset F_n$, the closure in
$\c$ of the endpoints of the axes of elements in
${\mathbb F}$ is denoted by
$C({\mathbb F})$ and is naturally identified with the
space of ends of ${\mathbb F}$.  Although the $UPG$ 
property referred to in the
following lemma is not defined until the next section, 
it is convenient to
place this result here.  This part of the lemma is  
quoted from
\cite{bfh:tits2} so there is no danger of circular 
reasoning.

\begin{lemma} \label{rank>1} Suppose that $\h$ is a
subgroup of Out$(F_n)$ and  that
${\mathbb F}$ is an $\h$-invariant (up to conjugacy)
subgroup of $F_n$ that is its own
normalizer.  Then
\begin{enumerate}
\item  There is a well-defined restriction
$\h|{\mathbb F}\subset$ Out(${\mathbb
F})$.
\item  If $\h$ is $UPG$, then $\h|{\mathbb F}$ is
$UPG$.
\item If ${\mathbb F}$ has rank at least two, then
any lift $\widetilde{\h|{\mathbb F}}
\subset Aut({\mathbb F})$  of $\h|{\mathbb F}$
extends uniquely to a lift
$\ti \h \subset Aut(F_n)$.
\item  If every element of a lift
$\widehat{\h|{\mathbb F}} \subset EH(C({\mathbb
F}))$
 of $\h|{\mathbb F}$  fixes at least three
points then $\widehat{\h|{\mathbb F}}$ extends
uniquely to a lift   $\hat \h \subset
EH(\c)$
\end{enumerate}
\end{lemma}

\noindent{\bf Proof of Lemma~\ref{rank>1}}  Since
${\mathbb F}$ is its own normalizer,
two automorphisms of $F_n$ that preserve
${\mathbb F}$ are conjugate by an element of
$F_n$ if and only if they are conjugate by an element
of ${\mathbb F}$. Part (1) follows
immediately.  Part (2) is Lemma 4.13 of
\cite{bfh:tits2}.  If ${\mathbb F}$ has rank at
least two, then there are no non-trivial inner
automorphisms that pointwise fix ${\mathbb
F}$.  It follows that the restriction homomorphism
from the subgroup of
$Aut(F_n)$ representing $\eta \in Out(F_n)$ to the
subgroup of
$Aut({\mathbb F})$ representing  $\eta|{\mathbb
F}$ is injective. Part (3) follows easily.
For part (4), note that   a non-trivial
covering translation does not fix more than two
points in $\c$ and so two lifts $\hat
\eta_1$ and $\hat \eta_2$ of $\eta$ that agree on
three points must be equal. \qed

\subsection{Relative train track maps}

    We study an element ${\cal \eta} \in Out(F_n)$
through its  lifts $\hat \eta
: \c \to \c$. We analyze the $\hat \eta$'s by
representing $\eta$ as a homotopy
equivalence $f : G \to G$ of a marked graph with
particularly nice properties
and studying the corresponding $\hat f$'s.  In this
subsection we recall some
properties of $f : G \to G$.

	  A {\it filtered graph} is a marked graph along
with a filtration  $\emptyset
=  G_0 \subset G_1 \subset \dots \subset G_K = G$
where each $G_i$ is obtained
from $G_{i-1}$ by adding a single edge $E_i$. We
reserve the words path and
loop for immersions of the interval and the circle
respectively. If $\rho$ is a
map of the interval or the circle into $G$ or
$\Gamma$, then $[\rho]$ is the unique
path or loop that is homotopic to $\rho$ rel
endpoints if any.  We say that a
homotopy equivalence
$f : G
\to G$ {\it respects the filtration} if each
$f(E_i) = E_iu_{i,f}$ for some loop $u_{i,f}  \subset
G_{i-1}$. Let $\he$\ be the group (Lemma 6.2 of
\cite{bfh:tits2}) of homotopy classes,
relative to vertices, of filtration respecting homotopy
equivalences of $G$.   There is
a natural map from $\he$\  to $Out(F_n)$.  We say
that $\eta \in Out(F_n)$ is $UPG$ (for
unipotent with polynomial growth) if it is in the
image of $\he$\  for some
$G$.  We say that a subgroup of $Out(F_n)$ is $UPG$
if each of its elements
is.  The main theorem of
\cite{bfh:tits2} states that every $UPG$ subgroup
$\H$  lifts to a subgroup
$\k$ of $FHE(G)$ for some
$G$.  We say that $\k$ is a {\it Kolchin
representative} of $\H$.

     In general we will use $\k$ to denote a subgroup
of $\he$.

         If
$\ti L
\subset
\Gamma$ is a line and $l$ is
the highest parameter value for which $\ti L$ crosses
a lift of $E_l$,
then   define the {\it highest edge splitting}
of $\ti L = \dots \tilde \sigma_{-1}\cdot \tilde
\sigma_0 \cdot \tilde
\sigma_1\dots$  by subdividing at the initial vertex
of each lift of $E_l$ that
is crossed, in either direction, by $\ti L$. We refer to 
the
vertices that
determine this decomposition as the {\it splitting
vertices} of the highest
edge splitting. If $\ti f : \Gamma \to \Gamma$ is a
lift of $f \in \he$ and if
$\ti L$ is
$\ti f_\#$-invariant, then Lemma 5.2
\ of
\cite{bfh:tits1}
implies that
$[\ti f(\ti
\sigma_j)]=
\ti \sigma_{j+r}$ for some $r$ and all $j$.  Roughly
speaking,
$\ti f$ acts on $\ti L$ by translating the highest edge
splitting by $r$ units.

   We will need the following  fixed point results.

\begin{lemma} \label{only one line}
\begin{itemize}
\item For any filtered graph $G$  and distinct 
$P_1,P_2,P_3 \in
\c$, there is a line $\ti L$ in $\Gamma$
connecting  $P_i$ to $P_j$ for some $1 \le i < j \le 3$
with the following property: If
$f
\in
\he$  and $P_1,P_2,P_3 \in Fix(\hat f)$ then $\ti f$
fixes each highest
edge splitting vertex in $\ti L$.
\item If $f \in \he$ and  $\ti f : \Gamma \to
\Gamma$ is  fixed point free,
then $\hat f$  fixes exactly two points.
\end{itemize}
\end{lemma}

\noindent{\bf Proof of Lemma~\ref{only one line}}
  For the first item, let  $\ti L_{i,j} \subset \Gamma$
be the line  connecting
$P_i$ to $P_j$ and  let $l_{i,j}$ be the highest
parameter value for which $\ti L_{i,j}$
crosses a lift of
$E_{l_{i,j}}$.  Assuming without loss that $l_{1,2} \ge
l_{1,3},l_{2,3}$, let $\ti L =
\ti L_{1,2}$ and let
$\ti v$ be any highest edge splitting vertex  of $\ti L$.
Suppose that $f \in \he$
and that $\hat f$ fixes each
$P_i$.  If $\ti f$ does not fix $\ti v$ then $\ti f$
translates the highest edge
splitting vertices  of
$\ti L$ away from one endpoint of $\ti L$, say $P_1$,
and toward the other,
$P_2$.  It follows that the highest edge splitting of
$\ti L$ is bi-infinite and hence
that
$l_1 = l_2 =l_3$.  This implies that $\ti f$ translates
the highest edge splitting
vertices of $\ti L_{2,3}$  away from $P_3$ and
toward $P_2$. But now $\ti f$ translates the highest
edge splitting vertices of $\ti
L_{1,3}$ away from $P_1$ and  away from $P_3$.
This
contradiction verifies the first item.

  We assume now that $f \in \he$ and that $\ti f$ is
fixed point free. By the first item, it suffices to show 
that $\hat f$ fixes
at least two points.
Implicit in Proposition 6.21 of \cite{bfh:tits1} is the
existence of  a half-infinite ray
$\ti R_+ \subset \Gamma$ with highest edge
splitting $\ti R_+ = \tilde \sigma_0
\cdot \tilde \sigma_1\dots$ such that $[\ti f(\ti
\sigma_j)]=
\ti \sigma_{j+r}$ for some $r>0$ and all $j \ge 0$.
Let $\ti v_j$ be the initial
vertex of $\ti \sigma_j$.   Since
$\ti f$ restricts to a bijection of vertices, there are
unique vertices
$\ti v_j$ so that $\ti f(\ti v_j) = \ti v_{j+r}$ for all
$j \in {\mathbb Z}$. Let
$\ti \sigma_j$ be the path connecting  $\ti v_{j}$ to
$\ti v_{j+1}$ for $j < 0$ and note
that
$[\ti f(\ti
\sigma_j)]= \ti \sigma_{j+r}$ for all $j$. If $l$ is
the highest parameter value for which $\ti R_+$
crosses a lift of $E_l$, then each $\ti
\sigma_j$ contains exactly one lift of $E_l$ and these
lifts are distinct for distinct
$j$. (The $j < 0$ case follows from the $j \ge 0$ case
which holds by
construction.) It follows that $\dots \tilde \sigma_{-
1}\cdot \tilde
\sigma_0 \cdot \tilde \sigma_1\dots$ is an
embedded line  whose
endpoints are both fixed by $\hat f$. \qed

\begin{lemma}\label{finite lifts}  Up to conjugation
by  covering translations, each
 $\eta \in Out(F_n)$ has only finitely
many lifts $\hat \eta : \c \to \c$ whose fixed point
set contains at least three points.
\end{lemma}

\noindent{\bf Proof of Lemma~\ref{finite lifts}}
Represent
$\eta$ by a homotopy equivalence $f : G \to G$ of
some marked graph and let $\ti f :
\Gamma \to
\Gamma$ be a lift whose extension $\hat f$ fixes at
least three points.  Given
$P_1,P_2,P_3
\in Fix(\hat f)$, let $\ti v \in \Gamma$ be the
unique vertex contained in
each of the  lines connecting  two of the $P_i$'s.  The
bounded cancellation lemma
\cite{co:bcc} implies that there is a bound, depending
on $f$ but not on the choice of
the lift $\ti f$ or the choice of points $P_i$, to the
length of the path $\ti \sigma$
connecting
$\ti v$ to
$\ti f(\ti v)$. In particular, the projected image
$\sigma$ takes on only finitely many
values as we vary the lift $\ti f$ and the choice of
points $P_i$.

			 Suppose that  $\ti f_i$, $i=1,2$,
are lifts of $f$ and that there
exist lifts $\ti v_i$ of a vertex $v$ and $\ti
\sigma_i$ of a path $\sigma$ such that
$\ti
\sigma_i$ is the path connecting $\ti v_i$ to $\ti
f_i(\ti v_i)$.  The covering
translation $T :
\Gamma \to \Gamma$  that carries $\ti v_1$ to
$\ti v_2$  satisfies $\ti f_2 T = T \ti
f_1$ and so conjugates $\hat f_1$ to $\hat f_2$.
\qed

\section{Property A}   \label{A}

     In this section we   prove that a finitely generated,
solvable
$UPG$ subgroup $\H$ embeds in  some $\Z^r$ and
hence that every solvable $UPG$
subgroup is free abelian.  Some of the arguments
proceed by induction on the
skeleta $G_i$ of a filtered graph $G$.   For this reason,
we do not assume in this
section that $\k$ is a  Kolchin representative of some
$\H$  but only that {\it $\k$ is
a finitely generated solvable subgroup of $\he$  with
the following feature: if $E_i$
is not a loop  then some $u_{i,f}$ is non-trivial.} In
particular, we do
not assume that  $\k$ injects into $Out(\pi_1(G))$
and we allow $G$ to have valence one
vertices. Note that if some $u_{i,f}$ is non-trivial,
then the terminal vertex of
$E_i$ is contained in a loop in $G_{i-1}$ and so must
have valence at least two in
$G_{i-1}$. Thus if $E_i$ is the first edge to contain a
vertex $v$, then either
$E_i$ is a loop or $E_i$ has $v$ as initial vertex.  In
either case, $v$
 is the initial vertex of $E_i$.

      The group of lifts $\ti f :
\Gamma \to \Gamma$ [respectively $\hat f : \c \to
\c$]  of elements of $\k$ has a
natural projection to $\k$.  By an {\it action of $\k$
on $\Gamma$ [respectively $\c$]
by lifts} we mean a section of this projection.  In
other words, an action by lifts  is an assignment $f
\mapsto \ti f$ [respectively $f
\mapsto \hat f]$ that respects composition.  Every
action $s$ of $\k$ on $\Gamma$ by
lifts determines an action $\hat s$ of $\k$ on $\c$ by
lifts and vice-versa.

    For each edge $E_i \subset G$ choose, once and for
all, a lift  $\tilde E_i^*
\subset \Gamma$.  Define $s_i(f) : \Gamma \to
\Gamma$ to be the  unique lift  of $f
\in \k$ that fixes the initial endpoint of $\ti E_i^*$
and note that {\it $s_i$  is an
action  of $\k$ on $\Gamma$ by lifts}.    Denote the
terminal endpoints of $E_i$ and
$\ti E_i^*$ by $v_i$ and $\ti v_i$ respectively.

\begin{defn} \label{propA}    If $E_i$ is not a 
component of $G_i$, denote the
component of
$G_{i-1}$ that contains $v_i$ by $B_i$ and the
(necessarily
$s_i(f)$-invariant) copy of the universal cover of
$B_i$ that contains $\ti
v_i$ by $\Gamma_{i-1}$.

      We say that $\k$ (and the choice of the $\ti
E_i^*$'s) satisfies {\it Property
A} (for abelian) if:
\begin{itemize}
\item For each $i$, either $E_i$  is a component of
$G_i$
 (in which case $E_i$ is a loop and $u_{i,f}$ is trivial
for all $f \in \k$)  or
$\ti v_i$ is a highest edge splitting vertex in a line
$\ti L_i \subset
\Gamma_{i-1}$ that is
$s_i(f)_\#$-invariant for all $f \in \k$.
\item If $\ti L_i$ and $\ti L_j$ have the same
projection in $G$,
then $\ti L_i =
\ti L_j$ and
$\ti v_i =
\ti v_j$.
\end{itemize}

\vspace{2in}

If $\ti L_i$ projects to an indivisible  loop
$\alpha_i$ and if some $u_{i,f}$ is non-trivial, then
we say (abusing notation
slightly) that $\alpha_i$ is an {\it essential axis} and
that $E_i$ is an {\it
essential edge}. The set of essential axes is denoted
${\it A(\k)}$ and the set of essential edges is denoted
${\it E(\k)}$.  If the
essential axis $\alpha$ is associated to $m_{\alpha}$
essential edges, then we
say that
$m_{\alpha}$ is the {\it multiplicity} of $\alpha$. In
the analogy with the mapping
class group of a compact surface, $A(\k)$ is the set of
reducing curves and $\k$ can be
chosen so that each $m_{\alpha} = 1$.

   For each essential axis $\alpha$, let $\ti
v_{\alpha}$ be the
preferred vertex of the preferred lift $\ti \alpha$
($=\ti v_i \in \ti L_i$
for any
$\ti L_i$ that projects to $\alpha$), let $T_{\alpha}$
be an indivisible
covering translation whose axis is $\ti \alpha$ and
let $s_{\alpha}(f)$  be
the lift of $f \in \k$ that fixes  $\ti v_{\alpha}$.
Note that {\it $s_{\alpha}$ defines an action of $\k$
on $\Gamma$ by lifts}.  Note also that if we think of
$\alpha$ as a closed path
with both endpoints at the projected image $v$ of
$\ti v$, then $[f(\alpha)] =
\alpha$.  It follows that $s_{\alpha}(f)_\#(\ti
\alpha) = \ti \alpha$ and hence
that $s_{\alpha}(f)$ commutes with $T_{\alpha}$.

\end{defn}

\begin{example}  \label{Dehn} We set notation for
this geometric example as shown
below: $M$  is the orientable genus two surface with
one boundary component; $A \subset
M$ is an embedded non-peripheral annulus  with 
boundary components $X$ and $E_1$;
$S = M \setminus A$; $G \subset M$ is
the embedded graph (spine of $S$) shown below with  
vertices $v_1$
and $v_2$  and  edges $E_1, \dots,  E_5$; $\H \cong
\mathbb Z$ is the subgroup
of the mapping class group of $M$ generated by the
Dehn twist  across the
annulus
$A$; and $\k \cong \mathbb Z \subset FHE(G,{\cal
V})$ is generated by $E_i
\mapsto E_i$ for $i \ne 2$ and $E_2 \mapsto
E_2E_1$.

\vspace*{2.5in}

       The universal cover of $M$ is identified with a
convex subset of the hyperbolic
plane
$H$ in the Poincare disk model.  Choose a lift $\ti L$
of $E_1$ and let $\ti A$ be the
lift of
$A$ that contains $\ti L$.  On
the left [respectively right] of $\ti A$ is a copy $\ti
S_1$ [respectively $\ti S_2$] of
the universal cover of $S$.  Given $D \in \H$, let
$\ti D_i : \ti M \to \ti M$  be the
lift that is the identity on $\ti S_i$..  Then $\hat
D_i$ fixes the endpoints of
$\ti S_i$  and 
$\hat D_1$ and $\hat D_2$ are the only lifts of $D$
that fix the endpoints of $\ti L$
and at least one other point. Note also that if $D$ is a
Dehn twist of order $k$ around
$A$ then $\ti D_1 = T^k \ti D_2$ where $T$ is the
indivisible covering translation that
preserves $\ti L$.

\vspace*{2.5in}

   On the graph level, there is a single essential edge 
$E_2$ with
essential axis $\alpha =
E_1$.  The lifts $\ti D_1$ and $\ti D_2$ correspond to 
$s_{\alpha}(f)$ and
$s_2(f)$ respectively.

   The free group $\pi_1(S)$ is 
generated by $E_2E_1 \bar E_2,
E_3E_1\bar E_3, E_4$ and $E_5$.   The  group $\H$
can be extended to a non-$UPG$
abelian group $\h$ by adding a generator that
restricts to a pseudo-Anosov homeomorphism
on $S$.  This can be represented by a relative train
map $f : G \to G$ that is the identity on $E_1$ and 
$E_2$ and that has a single
exponentially growing stratum with edges
$E_3,E_4$ and $E_5$.  \end{example}

\vspace{.1in}

   The following lemma justifies our notation.

\begin{lemma}\label{akr implies abelian}  If  $\k$
satisfies property A,
then there is an injective homomorphism
$\Phi_{\k} : \k \to {\mathbb Z}^r$
where
$r$ is the cardinality of $E(\k)$.  In particular, $\k$ is
free abelian.
\end{lemma}

\noindent{\bf Proof of Lemma~\ref{akr implies
abelian}}    If $E_i$ is an
essential edge with essential axis $\alpha$ then
$s_i(f)$ and $s_{\alpha}(f)$ are
lifts of $f$ that commute with $T_{\alpha}$ and so
differ by an iterate of
$T_{\alpha}$. Define a homomorpism $\phi_{\k}^i
: \k \to {\mathbb Z}$ by
$s_i(f) = T_{\alpha}^{\phi_{\k}^i(f)} s_{\alpha}(f)$
and define $\Phi_{\k}$ to be the
product of the
$\phi_{\k}^i$'s.   Then $f$ is in the kernel of
$\Phi_{\k}$ if and only
if $f(E_i) = E_i$ for each essential edge $E_i$ if and
only if $u_{i,f}$ is
trivial for each essential edge $E_i$.

  If $f \ne$ identity, then there is a smallest
parameter value $i>1$ for
which
$u_i(f)$ is non-trivial.		Since $f|G_{i-1}$ is
the identity, $s_i(f)$
restricts to a non-trivial covering translation of
$\Gamma_{i-1}$.  The line
$\ti L_i
\subset
\Gamma_{i-1}$ is $s_i(f)_\#$-invariant so must be 
the
axis of that covering
translation.  Thus $E_i \in E(\k)$  and
$f$ is not in the kernel of $\Phi_K$.  This proves
that $\Phi_K$ is
injective.
\qed

\vspace{.1in}

    The following lemma produces a pair of fixed
points in $\c$ or equivalently a fixed line in
$\Gamma$.

\begin{lemma} \label{two fixed points}  Suppose
that  $\hat \psi_1,
\dots, \hat \psi_m : C_\infty \to C_\infty$ are lifts
of elements of a
finitely generated UPG subgroup.  If the
$\hat
\psi_j$'s commute, then
$\cap_{j=1}^m$Fix($\hat \psi_j)$ contains at least
two points.
\end{lemma}

\noindent{\bf Proof of Lemma~\ref{two fixed
points}}  	     Choose a Kolchin
representative $\k$ of the UPG subgroup.   There are
elements
$f_j \in \k$ and commuting lifts $\ti f_j : \Gamma
\to \Gamma$ so that  $\hat f_j = \hat
\psi_j$ for $j=1,\dots,m$.

	    If $\ti f_1$ is fixed point free then
(Lemma~\ref{only one line}) Fix($\hat f_1)$
is a pair of points
$\{P,Q\}$.  Since
$\hat f_j$ commutes with $\hat f_1$, $\hat f_j$
setwise preserves
$\{P,Q\}$ and we need only show that $\ti f_j$ does
not reverse the
orientation on the line $\ti L$ connecting $P$ and
$Q$.  Let $\ti L = \dots \tilde
\sigma_{-1}\cdot \tilde \sigma_0 \cdot \tilde
\sigma_1\dots$ be highest edge splitting of $\ti L$. If
$\ti f_j$  reverses the
orientation on $\ti L$, then for some $j$, $[\ti f(\ti
\sigma_j)]$ equals
$\ti \sigma_j$ with its orientation reversed. The
projected image $\sigma_j$ determines
a conjugacy class in $F_n$ that is  periodic but not
fixed under the action of the outer
automorphism determined by $f$.   This contradicts
Proposition 4.5 of
\cite{bfh:tits2}.

    Suppose next  that  $Fix(\ti f_1) \ne \emptyset$
and that the group ${\cal T}(\ti
f_1)$ of covering translations that commute with
$\ti f_1$ is trivial.   The
fixed point set  of
$f_j$,    and hence of  $\ti f_j$, is a union of vertices
and edges. Since
${\cal T}(\ti f_1)$ is trivial,
each vertex
$v \in G$ has at most one lift $\ti v \in  Fix(\ti
f_1)$.  Since $\ti f_j$ commutes
with $\ti f_1$,  it preserves
$Fix(\ti f_1)$.  We conclude that each vertex in
$Fix(\ti f_1)$ is fixed by each
$\ti f_j$. (Recall that $f_j$ fixes each vertex in $G$.)
The  edges   with initial vertex in
$Fix(\ti f_1)$ project to distinct edges in $G$.  Let $\ti
E_k$ be the unique
such edge with minimal
$k$,  let $B'$ be the component of
$G_{k-1}$ that contains the terminal endpoint of
$E_k$ and let
$\Gamma_{k-1}' \subset \Gamma$ be the copy of
the universal cover of $B'$ that
contains the terminal endpoint of $\ti E_k$. Since
each $f_j$ maps an initial segment
of $E_k$ to an initial segment of $E_k$,
$\Gamma_{k-1}'$ is $\ti
f_j$-invariant for each $j$.  By our choice of $k$,  the
restriction $\ti
f_1|\Gamma_{k-1}'$ is fixed point free.   We can
now repeat the argument of the
first case on the restriction of the $\ti f_j$'s to
$\Gamma_{k-1}'$.

 \vspace{2in}

    Finally, suppose that  ${\cal T}(\ti f_1)$ is non-
trivial. Identify
${\cal T}(\ti f_1)$ with a subgroup ${\mathbb F}$ of
$F_n$.  The
space of ends of
${\mathbb F}$ is the closure
$C({\mathbb F}) \subset $ Fix($\hat f_1) \subset
C_{\infty}$ of the endpoints
of axes for elements of ${\cal T}(\ti f_1)$. Since $\ti
f_j$ commutes with
$\ti f_1$, the automorphism of ${\cal T}$
determined by $\ti f_j$ preserves
${\cal T}(\ti f_1)$ and
$\hat f_j = \hat \psi_j$ preserves $C({\mathbb F})$. 
By
Lemma~\ref{rank>1}, the
$\psi_j|{\mathbb F}$'s are contained in a  $UPG$
subgroup of Out(${\mathbb F}$).

      We argue by induction on
$m$, the $m=1$ case following from the fact that
$C({\mathbb F})$ contains at
least two points. Suppose  that $m > 1$.      By the
inductive hypothesis, there exist  $P,Q \in
C({\mathbb F})$ that are fixed by
the $m-1$ maps
$\hat f_2|C({\mathbb F}),\dots,\hat
f_m|C({\mathbb F})$.  Since $\hat
f_1|C({\mathbb F})$ is the identity,  $P$ and
$Q$ are fixed by each
$\hat f_j$.\qed

\vspace{.1in}

\begin{cor}\label{fg solv} Every finitely generated
solvable UPG subgroup $\H$  
 has a Kolchin representative that
satisfies property A; in particular,
$\H$ is  free abelian.
\end{cor}

\noindent{\bf Proof of Corollary~\ref{fg solv}}
Choose a Kolchin representative
$\k \subset \he$ for $\H$.     We work our way up 
the strata
$E_i$, modifying  $\k$ so that it satisfies
property A.  Denote the restriction of
$\k$ to $G_i$ by
$\k_i$.  Since $G_1$ is a single edge, $\k_1$ has
property A; we may assume by
induction  that $\k_{i-1}$ satisfies property A and is
therefore abelian.

  If $E_i$ is a
component of $G_i$ then $\k_i$ satisfies
property A. We may therefore assume that
$\Gamma_{i-1}$ and the action
$\hat s_i$ of
$\k$ on
$\c$ are defined.  The ends of $\Gamma_{i-1}$
define a subset $C_{\infty}^*
\subset C_{\infty}$.

 Let $R :
\k_i \to \k_{i-1}$ be the  restriction
homomorphism. If $R$ is an
isomorphism, then $\k_i$ is abelian.  Choose
generators
$g_1,\dots,g_m$ for $\k_i$ representing
$\psi_1,\dots,\psi_m \in \H$, and let
$\hat \psi_j^* =
\hat s_i(g_j)|C_{\infty}^*$. Note that the
$\hat \psi_j^*$'s are lifts of elements
of the $UPG$ subgroup $\H|\pi_1(G_{i-1})$ and
that, since $s_i$ is an action
and the $g_j$'s commute, the
$\hat \psi_j^*$'s  commute. Lemma~\ref{two fixed
points} therefore produces
$P_i,Q_i \in C_{\infty}^*$ that are fixed by each
$\hat s_i(g_j)$ and hence by
$\hat s_i(f)$ for each $f \in \k$.  The line  $\ti L_i$
connecting $P_i$ to $Q_i$
is $s_i(f)_\#$-invariant for each $f \in \k$.

     If $\ti L_i$ has the same image in $G$ as $\ti L_j$
for some $j < i$, then
after replacing $\ti E_i^*$ by a translate if necessary,
we may assume that $\ti
L_i = \ti L_j$.  Choose a highest edge splitting vertex
$\ti w_i$ of  $\ti L_i$
and apply the sliding operation of \cite{bfh:tits1}
simultaneously to
each
$f \in
\k$ reattaching $\ti E_i^*$ so that its terminal vertex 
is $\ti w_i$.   This does
not change
$s_i(f)|\Gamma_{i-1}$ and
respects the group structure of $\k$.  A new Kolchin
representative (still
called $\k$) is produced  that agrees with the old one
on $G_{i-1}$ and has the
additional feature that $\ti w_i$ is the terminal
vertex of $\ti E_i^*$; in
particular, the first condition of property A is
satisfied.  Since
$\ti w_i$ can be any highest edge splitting vertex, we
may assume without loss
that the second condition of property A is also
satisfied.

				  Suppose now that the
kernel $K$ of $R$ is non-trivial.  Each
$f^* \in K$ satisfies $E_j \mapsto E_j$ for $1 \le j
\le i-1$ and $E_i \mapsto
E_iu_{i,f^*}$.  Thus the restriction of $s_i(f^*)$ to
$\Gamma_{i-1}$
is  the covering translation determined by the lift of
$u_{i,f^*}$ beginning at
$\ti v_i$.  The assignment
$f^*
\mapsto u_{i,f^*}$ defines a  \lq suffix\rq\
homomorphism from $K$ into the free
group
$\pi_1(G,v_i)$.   Since
$\k$, and hence $K$, is solvable,  the image of this
homomorphism is
isomorphic to $Z$.  Thus the non-trivial
$s_i(f^*)|\Gamma_{i-1}$'s
have a common axis $\ti L_i$.  Choose $f^*$ so that
$s_i(f^*)|\Gamma_{i-1}$ is
non-trivial.        For each  $f \in
{\k}$, $s_i(ff^*f^{-1})|\Gamma_{i-1}$ is
a non-trivial covering translation with axis
$s_i(f)_\#(\ti L_i)$.
Since $K$ is normal in $\k$,
$s_i(f)_\#(\ti L_i) = \ti L_i$.  We have verified that
$\ti L_i$ is
$s_i(f)_\#$-invariant for each $f \in \k$.  The proof
now concludes as in
the previous case.
\qed

\section{Abelian subgroups are finitely generated}
\label{fin gen}

    In this section we prove that an abelian UPG
subgroup $\H$ is finitely
generated.  Section~\ref{A} and the exact sequence
$1\to\H\to\G\overset
\Omega \to\Z^b\to 1$ then implies that  every
solvable subgroup of Out($F_n$) is
finitely generated.  This fact was originally proved by
H. Bass and A. Lubotzky
\cite{bl:niltech}.  We also show that the above
sequence is a virtually
central extension.

    If $\k$ is a Kolchin representative of $\H$ that 
satisfies property A, $ f
\in
\k$ and $\Phi_{\k} : \k \to
{\mathbb Z}^r$ is the embedding of Lemma~\ref{akr
implies abelian}, then the
coordinates of  $\Phi_{\k}(f)$ are defined in terms of
the $s_i(f)$'s and the
$s_{\alpha}(f)$'s.  Our goal is to  recognize these lifts
of $f$ by their induced
actions on $C_{\infty}$, thus removing their
dependence on the choice of $\k$.

      We begin     placing an
additional restriction (the second item below) on our
Kolchin graphs and justifying our
assumption from section~\ref{A} (the first item
below).

\noindent \begin{lemma} \label{conditioned} For
every finitely generated UPG subgroup
$\H$ there is a Kolchin representative $\k$ with the
following properties.
\begin{itemize}
\item If $E_i$ is not a loop then  $u_{i,f}$ is non-
trivial for some $f \in \k$.
\item  Every vertex of $G$ is the initial vertex of at
least two edges.
\end{itemize}
\end{lemma}

\noindent {\bf Proof of Lemma~\ref{conditioned}}
Start with any Kolchin representative
$\k$.  After restricting to a subgraph if necessary, we
may assume that $G$ has
no valence one vertices.  We will show that if either
of the two properties fail, then we
can replace
$G$ by a graph with fewer edges.  This process
terminates after finitely many steps to
produce the desired Kolchin representative.

	    If the first property fails, then $E_i$ is
pointwise fixed  for each $f \in
\k$ and we may collapse it to a point.

						Suppose then
that the first property holds, that  $v$ is a vertex and
that
$E_i$ is the first edge that contains $v$.  If $E_i$ is
not a loop, then
 $u_{i,f}$ is a non-trivial loop  in
$G_{i-1}$  containing the terminal  endpoint of $E_i$
for some $f \in \k$.  Since loops
are assumed to be immersed, they cannot pass
through valence one vertices and $v$ must
be  the initial vertex of
$E_i$. For the same reason, $v$ must
also be the initial vertex of the second edge that is
attached to it.

    Suppose then that $E_i$ is a loop and that
$E_{j_0}, \dots ,E_{j_m}$
are the other edges that contain $v$.  If the second
property fails then
$E_{j_0}, \dots, E_{j_m}$ are non-loops with $v$ as
terminal endpoint and
$u_j(f) =  E_i^{k_j(f)}$ for $j=j_0,\dots,j_m$.
Redefine $\k$ by
replacing each $k_j(f)$ with $k_j(f)-k_{j_0}(f)$. This
can be achieved
on the graph level by sliding  $v$  around $E_i$
$j_0(f)$ times.  This
has no effect on the outer automorphism determined
by $f$ and respects the
group structure on $\k$.  The edge
$E_{j_0}$ is now fixed by each $f \in \k$ and so can
be collapsed to
point.	\qed

\begin{defn}  We say that a Kolchin representative
$\k$ for $\H$ is an
{\it \akr}  if it satisfies property A and the
conclusions of
Lemma~\ref{conditioned}.
\end{defn}

    We now turn to the task of determining, from the
action  of $\k$ on
$C_{\infty}$ by lifts,  if the axis of a given covering
translation projects  to an
element of $A(\k)$.

 Recall that in the analogy with the mapping class
group of a
compact surface $M$, $A(\k)$ corresponds to the set
of reducing curves in the minimal
reduction.  Such reducing curves are completely
characterized as follows  (See
Example~\ref{Dehn}, the proof of
Proposition~\ref{lifting to Aut} or
\cite{ht:surfaces}):  The free homotopy class
determined by a closed curve $\alpha
\subset M$ is an element of the minimal reducing
set for the mapping class represented
by a  homeomorphism $h : M
\to M$ if and only if for some (and hence each)
covering translation $T : \ti M \to \ti
M$ corresponding to $\alpha$, there are two lifts $\ti
h_1,
\ti h_2 : \ti M \to \ti M$ that commute with $T$
and whose extensions $\hat h_i$ over
the \lq circle at infinity\rq\ fix at least three points.
(If the free homotopy class
of $\alpha$ is fixed by $h$ but $\alpha$ is not one of
the reducing curves then there is
one such lift
$\ti h$.)   The analogous result in the non-geometric
case is given in
Corollary~\ref{independence of K}.

\begin{defn}    For any $\psi \in \H$ and covering
translation $T$,  define
$IL(\psi,T)$ (for Interesting Lifts) to be  the set of lifts
$\hat \psi : C_{\infty} \to C_{\infty}$ that commute
with $T$ and fix
 at least three points.
\end{defn}

    Recall that $\alpha \in A(\k)$ has a preferred lift
$\ti \alpha$ and an
indivisible covering translation $T_{\alpha}$ with
axis equal to $\ti \alpha$.
The next lemma states that the actions $s_{\alpha}$
and $s_i$ produce
interesting lifts.

\begin{lemma}\label{interesting}  Suppose that $\k$
is an \akr\ and that $\alpha$  is
the essential axis for $E_i \in E(\k)$.    Then
$\bigcap_{f \in \k}Fix(\hat
s_{\alpha}(f))$ and
$\bigcap_{f \in \k}Fix(\hat s_i(f))$ each contain at
least three points.  In
particular, if
$f$ represents $\psi \in \H$, then $\hat
s_{\alpha}(f),\hat s_i(f) \in
IL(\psi,T_{\alpha})$.
\end{lemma}

\noindent{\bf Proof of Lemma~\ref{interesting}}
There is a preferred
topmost splitting vertex $\ti v_{\alpha} \in\ti
\alpha$.
  Choose an edge $\ti E_j$ with initial vertex $\ti
v_{\alpha} $; by Lemma~\ref{conditioned}, we may
assume that $E_j \ne \alpha $.
If $u_{j,f}$ is   trivial for all $f \in \k$, then $E_j$ is a
loop that is fixed by
each
$f
\in \k$.   In this case, let $R_{\alpha}$  be an
endpoint of the axis for
$E_j$ that contains  $\ti E_j$.   If some
$u_{j,f}$ is non-trivial, then $\ti L_j$ is defined and
we
choose $R_{\alpha}$ to be an endpoint of the
translate $\ti L_j'$ of $\ti L_j$
associated to
$\ti E_j$.  In either case, each $\hat s_{\alpha}(f)$
fixes $R_{\alpha}$.

 \vspace*{3in}

  The proof for $s_i(f)$ is similar. By
Lemma~\ref{conditioned}
there exists an edge $\ti E_l \ne \ti E_i^*$, with the
same initial endpoint as
$\ti E_i^*$. Define $R_i \in$ Fix($\hat s_i(f))$ as in
the previous case using $ E_l$
in place of $ E_j$.\qed

\vspace{.1in}

    If $f \in \k$ represents $\psi \in \H$, then we use
$\hat s(f)$ and $\hat
s(\psi)$ interchangably.  We refer to the
$\hat s_i(\psi)$'s and the
$\hat s_{\alpha}(\psi)$'s as the {\it canonical lifts}
of $\psi$.  The next lemma and
corollary show that $A(\k)$ depends only on $\H$
and not on the choice of $\k$ and
that one can decide if  $\hat \psi$ is canonical from
its action on $\c$.

\begin{lemma} \label{axes are canonical}  Suppose
that $\k$ is an \akr\  and
that $\ti L$ is a line with endpoints $P,Q \in \c$.
Suppose further that:
\begin{itemize}
\item  $\ti f : \Gamma \to \Gamma$ is a lift of
some $f \in \k$
\item  Fix($\hat f$) contains $P,Q$ and at least one
other point.
\item  $\ti f$ does not fix the highest edge splitting
vertices of  $\ti
L$ .
\end{itemize}
Then there exists $i$ and a covering translation $T$
such that $T(\ti
L) =
\ti L_i$ and $\ti f = T^{-1} s_i(f) T $.
\end{lemma}

\noindent{\bf Proof of Lemma~\ref{axes are
canonical}}     By
Lemma~\ref{only one line}, Fix($\ti f) \ne
\emptyset$.  Choose an arc $\ti \sigma$
that intersects  Fix($\ti f$) only in its initial vertex,
say $\ti
p_1$, and that intersects the splitting vertices of $\ti
L$ only in its terminal
vertex, say $\ti v$.  Let $\ti E_i'$ be the first  edge of
$\ti
\sigma$ and let $\ti p_2$ be its terminal endpoint.
Since
$\ti f$ fixes
$\ti p_1$ but not
$\ti p_2$, $u_{i,f}$ is non-trivial and  $\ti L_i$ is
defined. Let $T$ be
the covering translation that carries $\ti E_i'$ to
$\ti E_i^*$ and let $\ti
L_i' = T^{-1}(\ti L_i)$ be the translate of $\ti L_i$
associated to $\ti E_i'$.
By construction, the ends  of $\ti L_i'$ are $\hat f$-
invariant.

\vspace{2in}

   If
$\ti L_i'
\ne
\ti L$ then there is an endpoint, say $S$, of $\ti L_i'$
that is
neither $P$ nor $Q$. Lemma~\ref{only one line}
implies that the highest
edge splitting vertices of either  $\ti L_{PS}$ (= the
line connecting $P$
to $S$) or
$\ti L_{SQ}$ are fixed. But  $\ti L_{PS}$  consists of a
segment of $\ti L$,
a segment of $\ti L_i'$ and perhaps a segment of
$(\ti \sigma \setminus \ti E_i)$. The last segment
contains no fixed
vertices.  By construction, the highest edge splitting
vertices of the
other segments are not fixed.  Thus the highest edge
splitting vertices of $\ti L_{PS}$ are not fixed.  The
symmetric  argument
applies to  $\ti L_{SQ}$ and yields the desired
contradiction. We conclude
that $\ti L_i' = \ti L$ and that $T(\ti v_i) = \ti p_2 =
\ti v$.  Since
$T^{-1} s_i(f) T$ and $\ti f$ both fix $\ti p_1$, they
must be equal.
    \qed

\vspace{.1in}

\begin{cor} \label{independence of K}If $\k$ is an
\akr\   of $\H$, then:
\begin{itemize}
\item  A covering translation $T$ corresponds to an
essential axis of $\k$ if and
only if
$IL(\psi,T)$ contains at least two elements for some
$\psi \in \H$.
\item For each $\alpha \in A(\k)$ and each $\psi
\in \H$,
$IL(\psi,T_{\alpha}) = \{\hat s_i(\psi): E_i$ is an
essential edge with axis
$\alpha\} \cup \{ \hat s_{\alpha}(\psi)\}$
\item $A(\k) = A(\H)$ depends only on $\H$ and
not on
$\k$.
\end{itemize}
\end{cor}

\noindent{\bf Proof of Corollary~\ref{independence
of K}} The
first and second items  are a direct consequence of
Lemma~\ref{interesting} and
Lemma~\ref{axes are canonical}. The third  item
follows from the first.
\qed

\begin{lemma} \label{fg} Every abelian UPG
subgroup ${\cal H}_0$ is finitely
generated.
\end{lemma}

\noindent{\bf Proof of Lemma~\ref{fg}}   Choose  an
increasing sequence
${\cal H}_1 \subset {\cal H}_2 \subset
\dots$ of finitely generated subgroups whose union
is  ${\cal H}_0$.  The
cardinality of $A(\h_j)$ and the multiplicities of the
elements of $A(\h_j)$ are
uniformly bounded (by the maximum number of
edges in a marked graph with the property
that the terminal vertex of each edge has valence at
least three.)  Since $A(\h_j)
\subset A(\h_{j+1})$, we may assume after passing
to a subsequence, that $A(\h_j)$
and the multiplicities are independent of $j$. In
particular, there is a fixed $r$ so
that
 $\Phi_j :{\cal H}_j \to {\mathbb Z}^r$ where
$\Phi_j$ is the embedding of
Lemma~\ref{akr implies abelian} and where ${\cal
H}_j$ has been identified with a
Kolchin representative $\k_j$.

    Fix $\psi \in \h_1, \alpha \in A(\h_1)$ and $j
\ge 1$; Let $IL(\psi,T_{\alpha})  =
\{b_1, b_2,\dots \}$.  For each $b_k$ and $b_l$, there
is an integer $p(k,l)$ such that
$b_lb_k^{-1} = T_{\alpha}^{p(k,l)}$.  For any fixed $f
\in \k$,
Corollary~\ref{independence of K} implies that each
coordinate
$\phi_j^i(f)$ of $\Phi_j(f)$ is one of the $p(k,l)$'s.
In
particular,
$\Phi_j(f)$ takes on only finitely many values as $j$
varies.

    Let $g_1,\dots,g_q$ be generators of ${\cal H}_1$.
After passing to a
subsequence we may assume that $\Phi_j(g_i)$ is
independent of $j$ for $i=1\dots
q$.  The lattice $\Phi_j({\cal H}_1)$ is therefore
independent of $j$.  It is
contained with finite index in a maximal lattice $L$
of rank $q$. The lattice
$\Phi_j({\cal H}_j)$ has rank $q$ and contains
$\Phi_j({\cal H}_1)$ so is contained in $L$.  In
particular, the index of
${\cal H}_1$ in
${\cal H}_j$ is uniformly bounded.  It follows that
${\cal H}_j = {\cal
H}_{j+1}$ for all sufficiently large $j$.

\qed

\vspace{.1in}

    We close this section by showing that
$1\to\H\to\G\overset
\Omega\to\Z^b\to 1$ is a virtually central
extension;  the proof is a variation on that
of Lemma~\ref{fg}.

\begin{lemma}\label{virtually central} $A(\H)$ is
$\h$ invariant (up to conjugacy).
There is  a finite index subgroup $\h' \subset \h$
whose actions on $\H$ and on $A(\H)$
are trivial.
\end{lemma}

\noindent{\bf Proof of Lemma~\ref{virtually
central}}  For each $\psi \in
\H$, each $\alpha \in A(\H)$ and each lift $\hat
\eta$ of each $\eta \in \h$,
conjugation by
$\hat \eta$ sends
$IL(\psi,T_{\alpha})$ to
$IL(\psi',T')$  where $\psi' = \eta \psi \eta^{-1}$ and
$T' = \hat \eta T_{\alpha}
\hat \eta ^{-1}$.  In particular, $IL(\psi,T_{\alpha})$
and $IL(\psi',T')$
have same cardinality.  Corollary~\ref{independence
of K} therefore
implies that the axis of $T'$ projects to an element of
$A(\H)$ and so
$A(\H)$ is invariant under the action of $\eta$.
Since
$A(\H)$ is finite, there is a finite index subgroup
for which this action is trivial.

  We  assume now that the action on $A(\H)$ is
trivial.  Choose $\hat
\eta$ so that
$T' = T_{\alpha}$ or equivalently so that $\hat \eta$
commutes with
$T_{\alpha}$.  Choose an \akr\
$\k$ for $\H$ and let $\Phi :\H \to {\mathbb
Z}^{r}$
be the embedding of Lemma~\ref{akr implies
abelian}.
 For any
pair of elements
$b_k,b_l
\in
IL(\psi,T_{\alpha})$, there  is an integer
$p(k,l)$ such that $b_lb_k^{-1} = T_{\alpha}^{p(k,l)}$.
Let $P(\psi)$ be
the finite collection of integers that occur as $p(k,l)$'s.
Since
conjugation by $\hat \eta$ carries
$IL(\psi,T_{\alpha})$ to
$IL(\psi',T_{\alpha})$  and since $\hat \eta $
commutes with $T_{\alpha}$,
$P(\psi)=P(\psi')$. By construction,
$\Phi(\psi')$ therefore takes on only finitely many
values as $\eta$ varies over
$\h$ and $\psi \in \H$ is fixed.  Since $\Phi$ is an
embedding, $\psi'$ takes on
only finitely many values.  After passing to a finite
index subgroup we may
assume that the action of
$\h$ by conjugation on $\psi$ is trivial.  After
applying this argument to a
finite generating set for $\H$, we see that the action
of $\h$ on $\H$
is virtually trivial.\qed

\section{Proof of Proposition~\ref{lifting to Aut}}
\label{lifting}

   The following lemma produces interesting lifts for
(iterates of)  individual elements
of $Out(F_n)$.

\begin{lemma} \label{lift one}  Suppose that $n \ge
2$ and that $\eta \in Out(F_n)$.
After replacing
$\eta$ by an iterate if necesssary,  there is  a lift
$\hat\eta \in EH(\c)$
   that fixes at least three points.  If $\gamma$ is a
non-trivial
primitive element of $F_n$ that is fixed (up to
conjugacy) by $\eta$ and if  $T$ is a
covering translation corresponding to  $\gamma$,
then we may choose $\hat\eta$ to commute
with $T$.
\end{lemma}

G. Levitt and M. Lustig inform us that they are 
developing techniques 
to understand the dynamics of
automorphisms of hyperbolic groups on the boundary 
of the group, and that
they are able to provide an alternate
proof of this lemma.

\noindent{\bf Proof of Lemma~\ref{lift one}}  We
assume at first that $\gamma$ and $T$
are given.

  The case that $\eta$ is
realized as an isotopy class of a surface
homeomorphism $h : S \to S$ is well
known (see for example Lemma 3.1 of
\cite{ht:surfaces}):  The Thurston classification
theorem implies, after replacing
$h$ by an iterate if necessary, that $S$ divides along
annuli into
 subsurfaces with negative Euler characteristic on
which $h$ is either the
identity or is pseudo-Anosov. Assuming that the
reduction is done along the
minimal number of annuli, we may choose the
curve representing $\gamma$ to lie
in one of the subsurfaces $S_i$.  If
$h|S_i$ is the identity, then it has a lift that fixes the
ends determined by
$S_i$ (See Example~\ref{Dehn}).  If $h|S_i$ is
pseudo-Anosov then $\gamma$
determines a boundary component of $S_i$ and
there is a lift fixing $\ti
\gamma$ and the endpoints of the singular leaves of
the pseudo-Anosov foliations
associated to that boundary component.

    We now turn to the general case and argue by
induction.  Since
every outer automorphism of $F_2$ is realized as a
surface isotopy class, the
preceding argument handles the $n=2$ case. We may
therefore assume that the
lemma holds for free groups of rank less than $n$.

    The smallest free factor ${\cal F}(\gamma)$ that
contains
$\gamma$ \cite{bfh:tits1} is $\eta$-invariant (up to
conjugacy). If $1 < $ rank$({\cal
F}(\gamma)) < n$, then the inductive hypothesis
provides a lift of $\eta|{\cal
F}(\gamma)$ with the desired properties.  Extending
this lift (Lemma~\ref{rank>1}) to
all of $F_n$ completes the proof.  We may therefore
assume that rank$({\cal F}(\gamma))$
is either 1 or
$n$.

      Choose an improved relative train track map  $f :
G \to G$
representing an iterate of $\eta$. Theorem 6.4
of \cite{bfh:tits1} contains a list of all the properties of
$f : G \to G$
that are used in this proof.
 The conjugacy class of $\gamma$ determines a loop
in $G$ that we also
call
$\gamma$.

    If $\cal F(\gamma)$ has rank $n$, then $\gamma$
must cross an edge in
the highest stratum of $G$.  If this stratum is
exponentially growing, then
Theorem 6.4 of  \cite{bfh:tits1}  implies that there is
an
$\eta$-invariant (up to conjugacy) subgroup
${\mathbb F}$ that is its own normalizer and
with the following additional property:  There is a
conjugacy between the the outer
automorphism  $\eta|{\mathbb F}$ and a
pseudo-Anosov mapping class $h : S \to S$.
Moreover, the conjugacy carries
$\gamma$ to a boundary component of $S$.  Since
${\mathbb F}$ has rank at
least two, lifts of $\eta|{\mathbb F}$ extend uniquely
to lifts of $\eta$. We
may therefore assume  that
${\mathbb F} = F_n$ and hence that
$\eta$ is represented by $h$. We are now reduced to
a previous case.

      If the top
stratum $G_m$ is not exponentially growing, then
$G_m$ is a single edge $E_m$.
The loop $\gamma$ splits (Lemma 5.2 of
\cite{bfh:tits1}) at the initial vertex $v$ of
$E_m$ each time that it crosses
$E_m$ in either direction.   We may therefore think
of $\gamma =\gamma_1\cdot
\dots \cdot \gamma_r$ as a concatenation of
Nielsen paths
 based at
$v$ (i.e. each $[f(\gamma_i)] = \gamma_i$).   We
claim that there are two distinct
$\gamma_l$'s. If not, then, since
$\gamma$ is indivisible, $\gamma = \gamma_1$ is
of the form $E_m \delta$, $\delta \bar
E_m$ or
$ E_m \delta \bar E_m$ for some path $\delta
\subset G_{m-1}$, where $\bar E_m$ is
$E_m$ with its orientation reversed.    In the first
two cases $\gamma$ is a free factor
and in the last case
$\gamma$ is freely homotopic to $\delta$.  Each of
these contradicts our assumption that
$\cal F(\gamma)$ has rank $n$ and so verify our
claim..

 The axis $Ax(T)  =
\dots
\ti
\gamma_1 \ti \gamma_2 \dots $ of $T$
 decomposes as a concatenation of lifts of the
$\gamma_l$'s. Let $\ti v$ be the initial vertex of
$\ti \gamma_1$ and let $\ti f$ be
the lift of $f$ that fixes $\ti v$.  Since  the
$\gamma_l$'s are  Nielsen paths, $\ti f$
fixes each concatenation point  in the decomposition
$Ax(T) = \dots \ti
\gamma_1 \ti \gamma_2 \dots $.  In particular,
$\hat \eta = \hat f$ fixes the
endpoints of
$Ax(T)$  and so commutes with $T$.  Now think of
$\gamma_1$ as a loop and extend $\ti
\gamma_1$ to the axis $Ax(T_1)$ of a covering
translation by concatenating
 translates of $\ti \gamma_1$.    Then $\hat \eta =
\hat f$ also fixes the endpoints
of $Ax(T_1)$.

   It remains to consider the case that  $\cal
F(\gamma)$ has rank one.  We may
assume that the first stratum $G_1 = \gamma$ is a
single edge; let $v$ be the
vertex of $G_1$.

      Suppose that an edge of an \egs\ $H_r$ is
attached to $v$. After replacing $f$ by
an iterate if necessary there is an edge $E$ of $H_r$
with initial vertex $E$ such
that $f(E) = E \cdot \beta$ splits into the
concatenation of $E$ with some non-trivial
path $\beta$.  Let
$\ti v$ be a lift of $v$ in the axis of  $T$, let  $\ti f :
\Gamma \to \Gamma$ be the
lift of $f$ that fixes $\ti v$ and let  $\ti E$ be the lift
of $E$ with initial vertex
$\ti v$.  Then
$\ti f(\ti E) =\ti E \cdot \ti \beta$ splits into the
concatenation of $\ti E$
and some other non-trivial path and so  each $[\ti
f^k (\ti E)]$  is a
proper initial segment of $[\ti f^{k+1} (\ti E)]$.  It
follows that  the
$[\ti f^k (\ti E)]$'s converge to an invariant ray
whose endpoint $R$ is fixed by
$\ti f$.  Let $\hat \eta = \hat f$.

       If there are no exponentially growing strata
attached to $v$, then there are no
zero strata attached to $v$ and each edge $E_i$ that
is  attached to $v$  is its own stratum $H_i$ and
satisfies $f(E_i) =E_iu_i$ for
some path $u_i \subset G_{i-1}$.  If $v$ is the initial
vertex of $E_i$ and
either $E_i$ is a loop or $u_i$ is non-trivial, then we
define
$\hat \eta$ as in the proof of
Lemma~\ref{interesting}:   Let $\ti v$ be
a lift of $v$ in the axis of  $T$, let  $\ti f : \Gamma
\to \Gamma$ be the
lift of $f$ that fixes $\ti v$ and let $\ti E$ be the lift of
$E$ with initial
vertex $\ti v$.  If $u_i$ is trivial then $R$ is the
endpoint of the axis for
$E_i$ that contains $\ti E_i$.  If $u_i$ is non-trivial
then $R$ is the
endpoint of the invariant ray $\ti u_i\cdot[\ti f(\ti
u_i)]\cdot[\ti f^2(\ti u_i)]\cdot
\dots$.

    In the remaining case, every $E_i$ attached to $v$
is either a fixed
non-loop or has $v$ as its terminal endpoint.
Proceeding as in the proof of
Lemma~\ref{conditioned} we can reduce the
number of edges in $G$ and arrive at
one of our previous cases.

     This completes the proof when $\gamma$ and
$T$ are given.  It remains to consider
the case that $\eta$ does not act periodically on any
conjugacy class in $F_n$.  By
induction on  $n$, we may assume that $\eta$ does
not act periodically on the
conjugacy class of any proper free factor in $F_n$.
We may therefore assume our improved
relative train track map $f : G \to G$ has only one
stratum and that this stratum is
exponentially growing.

 After passing to a further iterate if necessary, we may
assume
that there is a vertex $v$ and two edges $E_1$ and
$E_2$ initiating at $v$ such that
$f(E_i) = E_i
\cdot \beta_i$ splits into the concatenation of $E_i$
with some non-trivial path
$\beta_i$.  Choose a lift $\ti v$ of $v$, let $\ti f$ be
the lift of $f$ that fixes
$\ti v$ and let $\ti E_i$ be the lift of $E_i$ initiating
at $\ti v$.  As above,  $[\ti
f^k (\ti E_i)]$ converges to an invariant ray whose
endpoint $R_i$ is fixed by $\hat f$.
Moreover, the bounded cancellation lemma
\cite{co:bcc} and the fact that the lengths
of $[\ti f^k (\ti E_i)]$ grow exponentially in $k$
imply that $R_i$ is an  attracting
fixed point for the action of $\hat f$ on $\c$.

    We have shown that some iterate of $\eta$ has a
lift with at least two attracting
fixed points.  Applying this to $\eta^{-1}$, we
conclude (suppressing the iterate
for notational simplicity) that some $\hat{\eta}$ has
 at least two repelling fixed points.  If $\tilde\eta$ has
a fixed point, then
the preceding argument shows,after passing to an
iterate if necessary, that there are
also at least two attracting fixed points.  Suppose then
that $\tilde\eta$ is
fixed point free.  During the proof of Proposition 6.21
of \cite{bfh:tits1} we show
that there exists $\ti x \in \Gamma$ such that $\ti x,
\tilde\eta(\ti x),
\tilde\eta^2(\ti x) \dots$ is an infinite sequence in
an embedded ray $\ti B$;
the endpoint of $\ti B$ is fixed by $\hat\eta$ and is
not one of the
repelling fixed points. \qed

\vspace{.1in}

  We prove Proposition~\ref{lifting to Aut} in the
following equivalent form.

\vspace{.1in}

\noindent {\bf Proposition~\ref{lifting to Aut}} {\it
Every abelian subgroup
$\h \subset Out(F_n)$ has a virtual lift $\hat \h
\subset EH(F_n)$.  If $\gamma$ is
a non-trivial primitive element of $F_n$ that is
fixed, up to conjugacy, by each
element of $\h$ and if $T$ is a covering translation
corresponding to $\gamma$,  then
$\hat \h$ can be chosen so that each element
commutes with $T$}.

\vspace{.1in}

\noindent{\bf Proof of Proposition ~\ref{lifting to
Aut}}  We argue by induction on the
rank of the free abelian group $\h$.  If $\h$ has rank
one with generator $\eta$, then
$\hat \h$ is determined by choosing a lift $\hat
\eta$ that commutes with  $T$ if $T$ is
given.  We may now assume that the rank of $\h$ is
at least two and that the lemma holds
for all ranks less than that of $\h$.

    By Lemma~\ref{lift one}, there is an element $\eta
\in \h$ and  a lift $\hat \eta$
that fixes at least three points and that commutes
with $T$ if $T$ is given.  Let
$C = Fix(\hat \eta$) and let ${\cal T}(C)$ be the group
of covering
translations  that preserve $C$.  As described in
subsection~\ref{lifts},  $\hat \eta$
determines an automorphism $\Phi : F_n \to F_n$
whose fixed subgroup
${\mathbb F}$ corresponds to   ${\cal T}(C)$ under
the identification  of $F_n$ with
${\cal T}$. (See the proof of Lemma~\ref{T(C)} for
further details.) If
$\gamma$ and
$T$ are given, then $T \in {\cal T}(C)$ and
$\gamma \in {\mathbb F}$.

      We first show, after passing to a subgroup of $\h$
with finite index, that every
$\mu \in \h$ has a lift $\hat \mu$ that commutes
with $\hat \eta$ and with $T$ if $T$
is given.  In particular, $C$ is $\hat\mu$-invariant,
${\mathbb F}$ is $\h$-invariant
(up to conjugacy) and, if $\gamma$  is given then
each element of $\h|{\mathbb F}$ fixes
$\gamma$ up to conjugacy.

   Suppose at first that $\gamma$ and $T$ are not
given.  We say that two lifts
$\hat \eta_1$ and $\hat \eta_2$ of
$\eta$  are equivalent if $\hat \eta_1 = T_1 \hat
\eta_2 T_1^{-1}$ for some covering
translation $T_1$.  Let $IL(\eta)$ be the set of
equivalence classes of lifts $\hat \eta$
that fix at least three points. By  Lemma~\ref{finite
lifts}, $IL(\eta)$ is finite.
Since $\h$ is abelian, $\h$ acts by conjugation on
$IL(\eta)$.  After passing to a
subgroup of finite index, we may assume that this
action is trivial.  Thus for any lift
$\hat \mu_1$ of $\mu
\in \h$, there is a covering translation $T_1$ such
that $\hat \mu_1 \hat \eta \hat
\mu^{-1}_1 = T_1 \hat \eta T_1^{-1}$.  Thus  $\hat
\mu = T_1^{-1}\hat
\mu_1$  commutes with $\hat \eta$.

   Suppose now that $\gamma$ and $T$ are given.
Each $\mu \in \h$ has a lift $\hat
\mu$ that  commutes with $T$.  Since $\h$ is
abelian, $\hat
\mu \hat \eta \hat \mu^{-1}$ is a lift of $\eta$ that
commute with $T$ and so $\hat
\mu \hat \eta \hat \mu^{-1} = T^a\hat \eta$ for
some $a$.  Let $\hat \eta_k =
\hat \mu^k \hat \eta \hat \mu^{-k} =  \hat
\mu \hat \eta_{k-1} \hat \mu^{-1} = T^{ak}\hat
\eta$. Then each $\hat \eta_k$ is
conjugate to $\hat \eta$ and so fixes at least three
points. On the other hand, there
are only finitely many values of $l$ for which
$Fix(T^l\hat \eta) \ne Fix(T)$.  (This
follows from: (i) $\hat \eta$ fixes $Fix(T) = \{P,Q\}$
and so cannot move points very
near $P$ to points very near $Q$; and (ii)  $T$ acts co-
compactly on $C_{\infty}
\setminus  \{P,Q\}$.)  We conclude that $a = 0$ and
hence that $\hat \mu$ commutes
with $\hat \eta$.

   The proof now divides into cases, depending on the
rank of ${\cal T}(C)$.  If ${\cal
T}(C)$ is the trivial group, then  each $\mu$ has a
unique lift $\hat \mu$ that
preserves $C$ and so the assignment $\mu \mapsto
\hat \mu$ defines $\hat \h$.

    Suppose next that ${\cal T}(C)$ has rank one. Since
${\mathbb
F}$ is $\h$-invariant (up to conjugacy), we may
assume that $T$
is given and generates ${\cal T}(C)$.

     We claim that $C$ contains only finitely many
$T$-orbits.  This is a special case of the main theorem
of \cite{co:bcc}; the argument
in this case is short so we include it for completeness.
Let $P$ and
$Q$ be the endpoints of the axis of
$T$.  Since $T$ acts co-compactly on  $\c \setminus
\{P,Q\}$, it suffices to show that $P$ and $Q$ are the
only accumulation
points of $C$.  Suppose to the contrary that $S$ is an
accumulation
point other than $P$ and $Q$.  Arguing as in the
proof  of Corollary~\ref{finite
lifts}, using triples of points in $C$ limiting on $S$,
we find lifts $\ti v_i \in
\Gamma$ of a vertex $v \in G$ so that the covering
translation $T_i$ that carries
$\ti v_1$ to $\ti v_i$  commutes with $\hat \eta$.
But then $T_i$ is a
multiple of $T$ in contradictiction to the assumption
that $T_i(\ti
v_1) \to S$.  This verifies our claim.

    The action of $\h$ on the finitely many $T$-orbits
of $C$ is well
defined. After passing to a finite index subgroup, this
action is trivial. Choose a
point $R \in C \setminus\{P,Q\}$.  Composing
$\hat \mu$ with an iterate of $T$, we may  assume
that $\hat \mu$ fixes $R$.   The assignment $\mu
\mapsto \hat \mu$ defines
 $\hat \h$.

    Finally suppose  that ${\cal T}(C)$ has rank at least
two.  As noted above,
if $\gamma$ is given then it is fixed, up to
conjugacy, by each element of
$\h^* = \h|{\mathbb F}$.  Since
${\mathbb F}$ is the fixed subgroup of an
automorphism
$\Phi$ representing $\eta$, the image of
$\eta$ in
$\h^*$ is trivial.  Thus the rank of $\h^*$ is less
than that of $\h$ and by
induction, there is a lift $\widehat{\h^*} \subset
EH(C)$ to elements that commute with
$T$ if $T$ is given.  By Lemma~\ref{rank>1},
$\widehat{\h^*}$ extends to the
desired lift $\hat \h$.  \qed

\vspace{.1in}

\section{Proof of Theorem~\ref{main}}  \label{proof
of main}

     We may assume without loss that $\H$ is non-
trivial and, by
Lemma~\ref{virtually central},  that $\h$ acts
trivially on
$\H$ and on $A(\H)$. Let $\k$ be an
\akr\ for $\H$.  Throughout this section $\hat s =
\hat s_i$ or $\hat
s_{\alpha}$ and $T_{\alpha}$ is its
associated covering translation.

       Before proving the next lemma we show that it
implies Theorem~\ref{main}.

\begin{lemma} \label{extend} The lift $\hat s(\H)
\subset EH(\c)$  virtually extends to a
lift $\hat S(\h)\subset EH(\c)$ all of whose
elements  commute with
$T_{\alpha}$.
\end{lemma}

\vspace{.1in}

\noindent{\bf Proof of Theorem~\ref{main}} We
may assume that the extensions
$\hat S_i$ and
$\hat S_{\alpha}$   of $\hat s_i$ and $\hat
s_{\alpha}$ produced by
Lemma~\ref{extend} are defined on all of $\h$.
Define a homomorphism
$\phi_i' : {\cal H} \to \Z$ by $S_i(\psi) =
T_{\alpha}^{\phi_i'(\psi)}
S_{\alpha}(\psi)$ and note that the product $\Phi' :
{\cal H}
\to
\Z^r$ of the $\phi_i'$'s extends the embedding
$\Phi : \H \to \Z^r$ of
Lemma~\ref{akr implies abelian}.  Define $\Psi : \G
\to
\Z^{b+r}$ to be the product of $\Omega$ and
$\Phi'$.  Since $\H$ is the
kernel of $\Omega$ and $\Phi'|\H = \Phi$ is an
embedding,  $\Psi$ is an
embedding.

We now know that solvable subgroups $\G$ of
$\o(\f)$ are finitely generated and
 virtually abelian.
By passing to a subgroup of index at
most
$D(\n)$ where
$D(n):=|GL(\n,\Z/3\Z)|< 3^{n^2}$, we may assume
that the image of $\G$ in
 $\glthree$ is trivial.
Thus,
$\G$ has a subgroup of index at most this number
that is a torsion free
 Bieberbach group of $vcd$
at most
$vcd(\o(\f))=2\n-3$ (see
\cite{cv:moduli}). Also, a Bieberbach group of $vcd$
at most
$n$ has a subgroup of index at most $D(n)$ that is
free abelian (see, for
 example, \cite{lc:bg}).
Thus, $\G$ has a free abelian subgroup of index at
most $D(n)D(2n-3)< 3^{5n^2}$.
 This completes
the proof of the Theorem~\ref{main}.

\qed

\vspace{.1in}

   We are now reduced to Lemma~\ref{extend}.  The
proof uses Proposition~\ref{lifting to
Aut} and follows the general line of the proof of
Proposition~\ref{lifting to
Aut}.

  Let $C = \bigcap_{f \in \k}Fix(\hat s(f)$) and let
${\cal T}(C)$ be the group  of
covering translations that preserves $C$. Note that
${\cal T}(C)$ contains $T_{\alpha}$
and so has rank at least one. Let ${\mathbb F}$ be the
subgroup of $F_n$ corresponding
to ${\cal T}(C)$ and let $\ti s(\H) \subset Aut(F_n)$
be the lift of $\H$ corresponding
to $\hat s(\H) \subset EH(\c)$.

\begin{lemma}\label{T(C)}   ${\mathbb F}$ is the
fixed subgroup $\{ \gamma \in F_n : \ti
\eta(\gamma) =
\gamma$ for each $\ti \eta \in  \ti
s(\H)\}$ of $\ti
s(\H)$.
\end{lemma}

\noindent{\bf Proof of Lemma~\ref{T(C)}}  If the
endpoints of the axis
of  $T \in {\cal T}$ are contained in $C$, then they
are fixed by each $\hat
s(f)$ and so
$T$ commutes with each $s(f) : \Gamma \to
\Gamma$.  Thus $T$ commutes with each $\hat
s(f)$ and  each
$Fix(\hat s(f))$ is $T$-invariant. It follows that  $T
\in {\cal T}(C)$.

     Conversely, if $P$ and $Q$ are, respectively, the
backward  and forward endpoints
of the axis of $T$ then
$\lim_{n \to \infty}T^n(R) =Q$ and $\lim_{n \to
\infty}T^{-n}(R) =P$ for all $R \in \c
\setminus \{P,Q\}$.  If $T(C) = C$ then $C$ must
contain $P$ and $Q$.

    We have shown that $T \in {\cal T}(C)$ if and
only if the endpoints of the axis
of $T$ are contained in $C$.  By construction, the
latter condition is equivalent
to the endpoints of the axis of $T$ being fixed by each
$\hat s(f)$.  The lemma
now follows from the definition of $\ti
s(\H)$.\qed

\vspace{.1in}

\noindent{\bf Proof of Lemma~\ref{extend}} Given
$\mu \in \h$, choose a lift $\hat \mu$
that commutes with $T_{\alpha}$. Suppose also that
$\hat \eta \in \hat s
(\H)$ is given.  Since $\h$ acts trivially on $\H$,
$\hat
\mu
\hat
\eta
\hat \mu^{-1}$ is a lift of $\eta$ that commute with
$T_{\alpha}$ and so $\hat
\mu \hat \eta \hat \mu^{-1} = T_{\alpha}^a\hat
\eta$ for some $a$.  Arguing exactly as
in  the proof of Proposition~\ref{lifting to Aut}, we
conclude that $\hat \mu$ commutes
with $\hat \eta$. It follows that
$C$ is
$\hat
\mu$-invariant, that ${\mathbb F}$ is $\h$-
invariant up to conjugacy and  that each
element of
$\h|{\mathbb F}$ fixes $\alpha$ up to conjugacy.

   The proof now divides into cases, depending on the
rank of ${\cal T}(C)$. Suppose
that ${\cal T}(C)$ has rank one.     We claim that $C$
contains only finitely many
$T_{\alpha}$-orbits.   Let $P$ and $Q$ be the
endpoints of the axis of
$T_{\alpha}$.  Since $T_{\alpha}$ acts co-compactly
on  $\c \setminus
\{P,Q\}$, it suffices to show that $P$ and $Q$ are the
only accumulation
points of $C$.  Suppose to the contrary that $S$ is an
accumulation
point other than $P$ and $Q$.  Lemma~\ref{only one
line}, applied to triples of
points in $C$ limiting on $S$, implies that there are
vertices $\ti v_i \in \Gamma$
 that limit on $S$ and that are fixed by $\ti s(f)$ for 
each
$f \in \k$. There is no loss in
assuming that the $v_i$'s are  all lifts of the same
vertex in $G$.  The covering
translation
$T_i$ that carries
$\ti v_1$ to $\ti v_i$  commutes with each $s(f)$
and so must be a
multiple of $T_{\alpha}$. But this contradicts the
assumption that $T_i(\ti
v_1) \to S$.  This verifies our claim.

     There is a well defined action of $\h$ on the
finitely many $T$-orbits of $C$.
After passing to a finite index subgroup, this action is
trivial. By
Lemma~\ref{interesting} there exists $R \in \c$ that
is not an endpoint of the axis of
$T_{\alpha}$.  Composing
$\hat \mu$ with an iterate of $T$, we may  assume
that $\hat \mu$ fixes $R$.   The assignment $\eta
\mapsto \hat \eta$ defines
$\hat S$.

     We  may now assume that ${\cal T}(C)$ has rank
at least two.  By
Lemma~\ref{rank>1},  the restriction $\h^*$ of $\h$
to Out( ${\mathbb F}$) is well
defined. By Lemma~\ref{T(C)}, the restriction
$\H^*$ of $\H$ to Out(${\mathbb F}^*)$ is
trivial. Restricting the  exact sequence
$1\to\H\to\G\overset
\Omega \to\Z^b\to 1$ to ${\mathbb F}^*$ we see
that $\h^*$ is  abelian.     By
Propositon~\ref{lifting to Aut}, there is a virtual lift
$\widehat{S^*} \subset EH(C)$
of $\h^*$ such that each element of $\widehat{S^*}$
commutes with $T_{\alpha}$. Let $\hat
S$ be the unique extension of $\widehat{S^*}$ to a
virtual lift of $\h$.   Since
$\H^*$ is trivial,
$\hat S(\psi)$ restricts to the identity on $C$ for each
$\psi \in
\H$. Thus $\hat S(\psi)$ and
$\hat s(\psi)$ agree on   $C$ and so must agree
everywhere.  \qed


\providecommand{\bysame}{\leavevmode\hbox 
to3em{\hrulefill}\thinspace}

\end{document}